\documentclass{amsart}
\usepackage{amsmath,amssymb,amsfonts}
\usepackage{latexsym,amsthm}
\usepackage[dvips]{graphicx}
\usepackage{graphicx}
\usepackage[usenames]{color}
\usepackage{rotating}
\usepackage{hyperref}
\usepackage{url}

\begin{document}
\title{FRAY JUAN DE ORTEGA's APPROXIMATIONS, 500 years after}


\author[M. Benito]{Manuel Benito
\address{ Logro{\~n}o, La Rioja, Spain}
\email{mbenit8@palmera.pntic.mec.es}}

\author[J. J. Escribano]{Jose Javier Escribano
\address{Calahorra, La Rioja, Spain}
\email{jjescribanob@gmail.com}}

\author[E. Fern\'andez]{Emilio Fern\'andez
\address{UR, Logro{\~n}o, La Rioja, Spain}
\email{eferna35@gmail.com}}

\author[M. S\'anchez]{Mercedes S\'anchez
\address{UCM, Madrid, Spain}
\email{merche@mat.ucm.es}}



\renewcommand{\abstractname}{Abstract}
\begin{abstract}
 In 1512, on December $30^{th}$, the first edition of  Fray Juan de Ortega's Arithmetic was published in Lyon. The last chapter, titled ``Rules of Geometry", deals with lower approximations of 14 square roots.
 In later editions of the Arithmetic on 1534, 1537 and 1542 in Seville,
 these values are replaced  by upper approximations.
 Twelve of them verify the Pell's equation, and so they are optimal.
 At this moment nobody knows the way they were obtained. In this paper we show how these approximations can be obtained through a method consistent with the mathematical knowledge at that time.

\vspace{1cm}

{\noindent\scshape R\'esum\'e.\ }
Fray Juan de Ortega
publia \`{a} lyon  la premi\`ere \'{e}dition de son
 Arithm\'etique le 30 D\'ecembre 1512.
 Dans son dernier chapitre intitul\'e ``R\`egles de la g\'eom\'etrie" il y a des approximations par d\'efaut de 14 racines carr\'ees. Dans les \'editions suivantes de S\'eville de 1534, 1537 et 1542, ces valeurs ont \'et\'e remplac\'ees par des approximations par exc\`es. Douze d'entre elles sont optimales, c'est \`a dire v\'erifient l'\'equation de Pell.
 M\^eme aujourd'hui, personne ne sait encore comment elles ont \'et\'e obtenues. Cet article d\'ecrit une m\'ethode pour obtenir  ces approximations, compatible avec les connaissances math\'ematiques de l'\'epoque.

\vspace{1cm}

{\noindent\scshape Resumen.\ }
El 30 de diciembre de 1512, Fray Juan de Ortega  public\'o en Lyon la primera edici\'on de su Aritm\'etica. En el \'ultimo cap\'itulo, ``Reglas de geometr\'ia", aparecen aproximaciones por defecto de 14 ra\'ices cuadradas. En las ediciones de Sevilla de 1534, 1537 y 1542, se sustituyen estos valores por aproximaciones por exceso. Doce de ellas son \'optimas (verifican la ecuaci\'on de Pell). Hasta la fecha se desconoce como fueron obtenidas. En este art\'iculo, se expone un m\'etodo por el que se obtienen todas estas aproximaciones y es coherente con los conocimientos matem\'aticos de la \'epoca.

\end{abstract}

\maketitle


 \section{Introduction}
  \begin{center}
\begin{figure}
\includegraphics[scale=0.82]{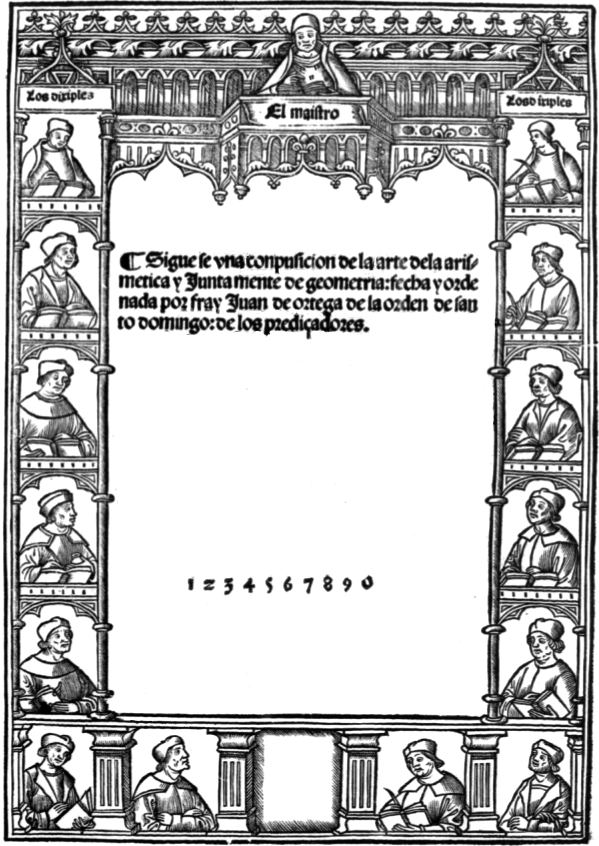}
\caption{Front cover. 1512 edition.}
\end{figure}
\begin{figure}
\includegraphics[scale=0.7]{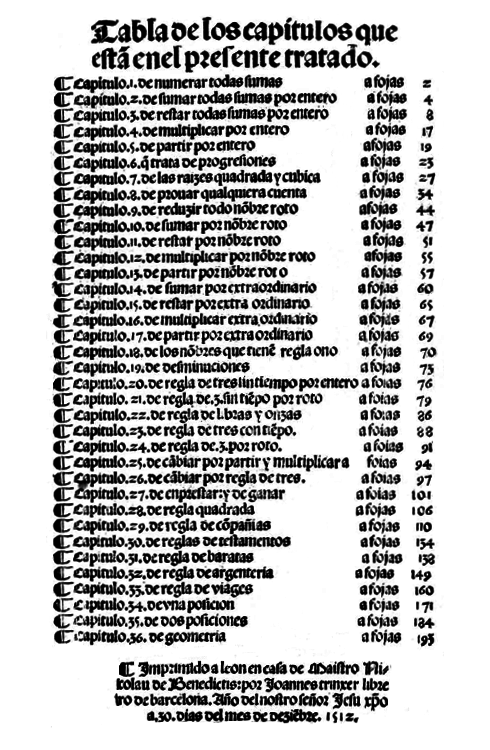}
\caption{Table of Contents.  1512 edition.}
\end{figure}
\end{center}

\begin{figure}
\includegraphics[scale=0.5]{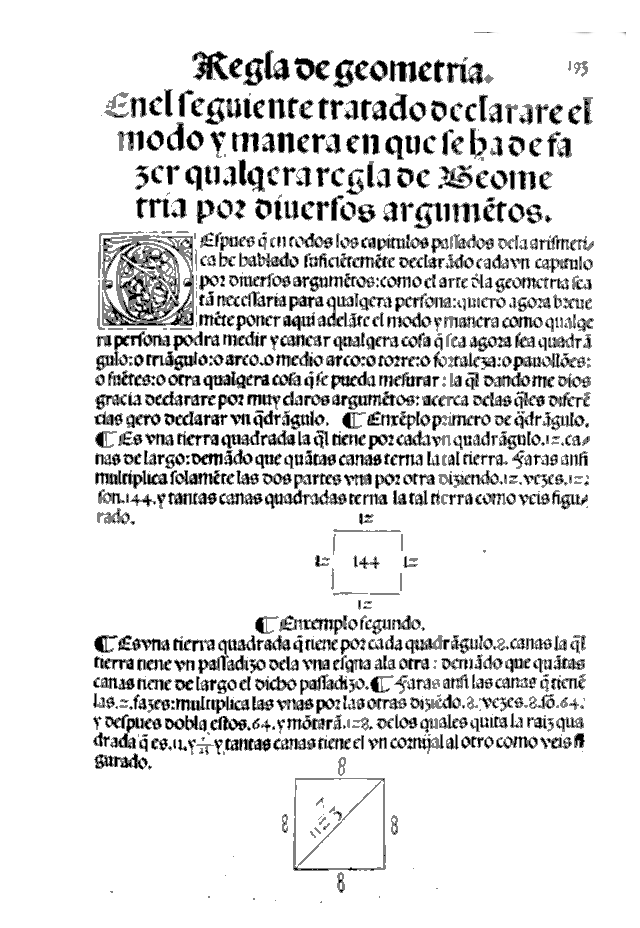}
\caption{Page 193. 1512 edition.}
\end{figure}

 The first edition of Fray Juan de Ortega's Arithmetic\footnote{ He was born in Palencia about 1480. He belonged to the religious order of preachers and was sent to the province of Arag\'on.

 He devoted his life to the teaching of Mathematics in Spain and Italy.} was published on December $30^{th}$ 1512 in Lyon. The book contains a collection of elemental arithmetic rules, such as operations with integers and rational numbers, square roots,\ldots, some notions  about commercial calculus such as proportions or equivalences between Spanish coins of that time.
  The rules are written in a
  practical and didactic way, so ``there will not be any fraud in the world about computing".

The text achieved a great success in Europe. In 1515 the work was published and translated into Italian in Roma and it was also translated into French by Claude Plantin and published in Lyon.  After this edition there were some new ones editions: in Messina (1522), in Seville (1534, 1537, 1542, 1552), in Grenade (1563). However, we remember this book, 500 years later, because of the approximations of the square roots written on the geometric applications at the end. In the last chapter, ``Rules of Geometry", Ortega solved exercises of elementary geometry and had to find some square roots. Some of them, in our language, are about  findding the length of a field with a circular shape and whose area is equivalent to another one with a square shape, or to find the edge of an equilateral triangle so that its area is equivalent to a given square. Both problems are written in the Appendix II.

In the first edition he approaches 14 lower square roots following a rule previously exposed.  On the Seville editions in 1534, 1537 and 1542, those values were replaced, without any explanation, by upper approximations optimal in 12 cases out\footnote{The upper approximation  $\frac xy$ of $\sqrt{n}$ is optimal if $(\frac xy)^2-n=\frac {1}{y^2}$, these equations are known as Pell'equations. A well known example of this equation is the problem of Archimedes'cows. In the XII century, some similar equations were solved by Bascara.}.

In the edition of 1512, page 230, there are also two upper approximations $\sqrt{127\frac{3}{11}}\simeq 11\frac{2}{7}$ y $\sqrt{5\frac 13}\simeq 2+\frac 16+\frac 17$, that they haven't been changed in the following editions\footnote{See Appendix II}.

How did Ortega manage to obtain these values? It has been a mystery that has occupied  the mind of many mathematicians and historians of science and has led to a lot of papers. Some of them are due to P. Tannery \cite{tan1}, \cite{tan2},\cite{tan3},\cite{tan},  Cantor \cite{can}, Enestr\H{o}m \cite{ene}, Perrot \cite{pert}, J. Rey Pastor \cite{rey1}, \cite{rey2}, \cite{rey3}, \cite{rey4}, \cite{rey5}, and Barinaga \cite{bar1}, \cite{bar2}.

Some of the results presented by Ortega  could not be obtained by the methode used by those authors to calculate the roots,  (they considered them mistakes of the author or the printer).\footnote{They also fail to discuss the upper approximations $\sqrt{5\frac 13}\simeq 2+\frac 16+\frac 17$ and $\sqrt{127\frac{3}{11}}\simeq 11\frac{2}{7}$ which remain unmodified in all editions and that can be seen in Appendix II.}

Although the solutions obtained by Ortega in the three editions of Seville are optimal, they should be something very new to the point to change the values; when Gonzalo del Busto re-edit the work, he was forced to rectify the many mistakes he found in some previous impressions. Thus, the optimal approximations were replaced by those from the 1512 edition.
In the following pages we will show all the roots and the approximations proposed by Ortega in different editions. And, given to the absence of explanations by the author, we suggest what could happen and we present a method to obtain all approximations, consistent with the mathematical knowledge of the time.

\section{Relationship between the roots and the approximations proposed by Ortega}

The roots that appear in the work are the following:
$$\begin{array}{r|c|c|c|c|c|}
\hline
& & \text{I}& &\text{II}& \\
\hline
 & \text{Number}& \text{Editions}\ldots&\text{Rest}&\text{ Editions}&\text{Rest}\\
 && 1512 \ldots&&1534-1537-1542&\\
 \hline
 1&128&11\frac{7}{23}&\frac{112}{529}&11\frac{16}{51}&\frac{-1}{2601}\\
 \hline
 2&80&8\frac{16}{17}&\frac{16}{289}&8\frac{17}{18}& \frac{-1}{324}\\
 \hline
 3&297&17\frac{8}{35}&\frac{216}{1225}&17\frac{659}{2820}&\frac{-1}{7952400}\\
 \hline
 4&300&17\frac{11}{35}&\frac{264}{1225}&17\frac{25}{78}&\frac{-1}{6084}\\
 \hline
 5&375&19\frac{14}{39}&\frac{350}{1521}&19\frac{285}{781}&\frac{-1}{609961}\\
 \hline
 6&135&11\frac{14}{23}&\frac{126}{529}&11\frac{13}{21}&\frac{-1}{441}
\\
 \hline
 7&75&8\frac{11}{17}&\frac{66}{289
}&8\frac{103}{156}&\frac{-1}{24336}\\
 \hline
 8&756&27\frac{27}{55}&\frac{ 756}{3025
}&27\frac{109}{220}&\frac{-1}{48400}\\
 \hline
 9&611&24\frac{35}{49}&\frac{10}{49}&24\frac{6886}{9585}&\frac{-1}{91872225}\\
 \hline
 10&231&15\frac{6}{31}&\frac{150}{961
}&15\frac{151}{760}&\frac{-1}{577600}\\
 \hline
 11&800&28\frac{16}{57}&\frac{656}{3249}&28\frac{197}{693}&\frac{-1}{480249}\\
 \hline
 12&4100&64\frac{4}{129}&\frac{ 500}{16641}&64\frac{1}{32}&\frac{ -1}{1024
}\\
 \hline
 13&2000&44\frac{64}{89}&\frac{1600}{7921}&44\frac{2079}{2882}&\frac{-89}{68644
}\\
 \hline
 14&9600&97\frac{191}{195}&\frac{764}{38025}&97\frac{191}{194}&\frac{-36481}{37636}\\
 \hline
 15&127\frac{3}{11}&11\frac{2}{7}&\frac{-51}{539}&11\frac{2}{7}&\frac{-51}{539}\\
 \hline
 16&5\frac{1}{3}&2+\frac{1}{6}+\frac{1}{7}&\frac{-1}{1764}&2+\frac{1}{6}+\frac{1}{7}&\frac{-1}{1764}\\
 \hline
 \end{array}
$$

In columns, the approximations of the first edition of 1512, french traslation of 1515, and the approximations modified in following editions are written

All first 14 values are lower approximations. They have been obtained by applying the usual algorithm for square roots described in chapter 7 of the book titled ``About square and cubic roots". In our current language
$$a+\frac{r}{2a+1} \leq \sqrt{n}.$$

Where  $a=\lfloor \sqrt{n}\rfloor$ (integer part) and $r=n-a^2$, the rest.

Thus, $$\sqrt{128}\simeq 11+\frac{128-121}{22+1}=11+\frac{7}{23}.$$

The last two values in column I are upper approximations and the way they were obtained is unknown. Both remain unchanged in later editions.

In column II all values that appeared in Seville's editions in 1534, 1537 y 1542 are written. They are upper approximations and all of them are optimal except those that are written in rows 13, 14 and 15. In those cases we also indicate how the approximations can be done.

The roots that appear in the Italian editions are
 $$
\begin{array}{|r|c|c|c|}
\hline
 & \text{Number}& \text{Roma 1515}&\text{Rest}\\
 \hline
 1&108&10+\frac{8}{21}&\frac{104}{441}\\
 2&135&11\frac{14}{23}&\frac{126}{529}\\
 3&187&13+\frac{2}{3}&\frac{2}{9}\\
 4&756&27\frac{27}{55}&\frac{ 756}{3025}\\
 5&611&24\frac{35}{49}&\frac{10}{49}\\
 6&800&28\frac{16}{57}&\frac{656}{3249}\\
 7&75&8\frac{11}{17}&\frac{66}{289}\\
 8&128&11\frac{7}{23}&\frac{112}{529}\\
 9&127+\frac{3}{11}&11\frac{2}{7}&\frac{-51}{539}\\
 10&231&15\frac{6}{31}&\frac{150}{961}\\
 11&5+\frac{1}{3}&2+\frac{1}{6}+\frac{1}{7}&\frac{-1}{1764}\\
 \hline
 \end{array}
 \
\begin{array}{|r|c|c|c|}
\hline
 & \text{Number}& \text{Mesina 1522}&\text{Rest}\\
 \hline
 1&75&8+\frac{11}{17}& \frac{66}{289}\\
 2&91&9+\frac{10}{19}&\frac{90}{361}\\
 3&871&29+\frac{30}{59}&\frac{870}{3481}\\
 4&4644&68+\frac{20}{137}&\frac{2340}{18769}\\
 5&40+\frac{1}{2}&6+\frac{9}{26}&\frac{153}{676}\\
 6&33+\frac{1}{3}&5+\frac{25}{33}&\frac{200}{1089}\\
 7&20+\frac{4}{11}&4+\frac{16}{33}&\frac{2072}{1089}\\
 8&18+\frac{3}{4}&4+\frac{11}{36}&\frac{275}{1296}\\
 \hline
 \end{array}\label{em}
 $$

\section{Our hypothesis about the method  applied by Ortega}

Rey Pastor in \cite{rey2} page 80, point out the possibility that Ortega was inspired on Nicolas Chuquet' \textit{Triparty}, or on any book of Arabic origin, even on the Paciolo' \textit{Summa}.

Our hypothesis is that Ortega could obtained his approaches using a special method of ``la regle des nombres mohines"(``mediation"\footnote{A mediation of two fractions $\frac{a}{b}$ and $\frac{c}{d}$, $a,b,c,d >0$, is the  fraction $\frac{a+c}{b+d}$ between both fractions.} , some kind of mean) written in the manuscript \textit {Triparty en la sciencia des nombres}, from  Chuquet \cite{chu}. This book was finished in 1484\footnote{ However the text was not  printed until 1880 by Aristide Marre. See Chuquet \cite{chu1}} in Lyon, the city where Ortega published the first edition of his work in 1512.

The Chuquet's text contains the approximations of the square roots of the first natural numbers. Chuquet began with the integer part of the lower and upper square root and apply in a systematical way his rule. Here the approximations of square roots of 14 first natural numbers are computed\footnote{This method provides solutions that verify the Pell's equation.}.

For $\sqrt{6}$ he obtains the following approximations
$$2+\frac{0}{1}<\sqrt{6}<2+\frac{1}{1}$$
$$2+\frac{0}{1}<\sqrt{6}<2+\frac{1}{2}$$
$$2+\frac{1}{3}<\sqrt{6}<2+\frac{1}{2}$$

$$\cdots $$

 $$\begin{array}{|c|c|c|}
 \hline \\
 \text{Number}&\text{Approximation
of the square root}& \text{Upper error}\\
 \hline
 2&1+\frac{169}{408}& \frac{1}{166464}\\
 3&1+\frac{571}{780}& \frac{1}{608400}\\
 5 & {\color{green}2+\frac{161}{682} }& \frac{1}{465124}\\
  6 & 2 +\frac{881}{1960} & \frac{1}{3841600}\\
7  & 2+\frac{7873}{12192} & \frac{1}{148644864}\\
8  & 2+\frac{985}{1189} & \frac{1}{1413721}\\
10  & 3+\frac{1405}{8658} & \frac{1}{74960964}\\
10  & 3+\frac{228}{1405} & \frac{-1}{1974025}\\
11  &3 +\frac{379}{1197} & \frac{1}{1432809}\\
12  & 3+\frac{181}{390} & \frac{1}{1432809}\\
13  & 3+\frac{109}{180} & \frac{1}{32400}\\
 14  &3 +\frac{2667}{3596} & \frac{1}{12931216}\\
 \hline
 \end{array}
 $$

By applying the Chuquet's method to 16 values in column II, we reach  Ortega's solution in 14 cases. However by this way we can't obtain the solution for $\sqrt{9600}$ and, for $\sqrt{2000}$, the value obtained is $44\frac{189}{262}$, which is the result  of simplifying by 11 the Ortega's solution ${\color{red}44 \frac{2079}{2882}}$.

We suppose that Ortega could have used the mediation rule taking as lower initial value the value given in the 1512 edition and the upper value which is obtained by decreasing the denominator in one unit\footnote{Nowadays the inequality is well known, however the second part of the inequality (Heron's formula) was not known until the XIX century.See \cite{tan1}} as we can see below. In our language
$$a+\frac{r}{2a+1} \leq \sqrt{n}\leq a+\frac{r}{2a}. $$
As above, $a=\lfloor \sqrt{n}\rfloor$ and, $r=n-a^2$ the lower rest.
Thus for the first value of $\sqrt{128}$ the following approximations are obtained:
$$11+\frac{7}{23}<\sqrt{128}<11+\frac{7}{22}$$
$$11+\frac{14}{45}<\sqrt{128}<11+\frac{7}{22}$$
$$11+\frac{21}{67}<\sqrt{128}<11+\frac{7}{22}$$
$$\cdots $$

In the following tables we show the results obtained by this way for the 14 first values of column II. The red values are the solutions given by Ortega, even though we have continued the process to find the first optimal solution.

As we have seen, the method provides all of Ortega's solutions, even the appro \- ximation of $\sqrt{2000}$ which again appears simplified. Later we will see another way to approach $\sqrt{2000}$ as Ortega did (without simplifying).

The approximations written in the two last rows $\sqrt{127\frac{3}{11}}$ and $\sqrt{5\frac{1}{3}}$ have three singular things in relation to those given above.

Those are upper approximations that remain unchanged in all editions and are not integers. As the edition of 1512 does not gave a lower value to start the process described above, we have started from the lower and upper integer roots as Chuquet did.

 $$\begin{array}{|c|c|c|c|}
 \hline
 {\bf 128}&\text{Lower approximation}&\text{Upper approximation}&\text{Upper error}\\
 \hline
1&11+\frac{7}{23}&11+\frac{7}{22}&\frac{49}{484}\\
2&11+\frac{14}{45}&11+\frac{7}{22}&\frac{49}{484}\\
3&11+\frac{21}{67}&11+\frac{7}{22}&\frac{49}{484}\\
4&11+\frac{21}{67}&11+\frac{28}{89}&\frac{161}{7921}\\
5&11+\frac{21}{67}&11+\frac{49}{156}&\frac{217}{24336}\\
6&11+\frac{21}{67}&11+\frac{70}{223}&\frac{217}{49729}\\
7&11+\frac{21}{67}&11+\frac{91}{290}&\frac{161}{84100}\\
8&11+\frac{21}{67}&{\color{red}11+\frac{16}{51}}&\frac{1}{2601}\\
\hline
\end{array}$$

 $$\begin{array}{|c|c|c|c|}
 \hline
 {\bf 80}&\text{Lower approximation}&\text{Upper approximation}&\text{Upper error}\\
 \hline
 1&8+\frac{16}{17}&9+\frac{0}{1}&\frac{1}{1}\\
2&8+\frac{16}{17}&{\color{red}8+\frac{17}{18}}&\frac{1}{324}\\
 \hline
\end{array}$$

$$\begin{array}{|c|c|c|c|}
 \hline
 {\bf 297}&\text{aproximaci\'on por defecto}&\text{aproximaci\'on por exceso}&\text{error por exceso}\\
 \hline
 1&17+\frac{8}{35}&17+\frac{4}{17}&\frac{16}{289}\\
2&17+\frac{3}{13}&17+\frac{4}{17}&\frac{16}{289}\\
3&17+\frac{7}{30}&17+\frac{4}{17}&\frac{16}{289}\\
4&17+\frac{7}{30}&17+\frac{11}{47}&\frac{27}{2209}\\
5&17+\frac{7}{30}&17+\frac{18}{77}&\frac{16}{5929}\\
6&17+\frac{25}{107}&17+\frac{18}{77}&\frac{16}{5929}\\
7&17+\frac{25}{107}&17+\frac{43}{184}&\frac{9}{33856}\\
8&17+\frac{68}{291}&17+\frac{43}{184}&\frac{9}{33856}\\
9&17+\frac{111}{475}&17+\frac{43}{184}&\frac{9}{33856}\\
10&17+\frac{154}{659}&17+\frac{43}{184}&\frac{9}{33856}\\
11&17+\frac{154}{659}&17+\frac{197}{843}&\frac{31}{710649}\\
12&17+\frac{154}{659}&17+\frac{351}{1502}&\frac{37}{2256004}\\
13&17+\frac{154}{659}&17+\frac{505}{2161}&\frac{27}{4669921}\\
14&17+\frac{154}{659}&{\color{red}17+\frac{659}{2820}}&\frac{1}{7952400}\\
 \hline
\end{array}$$

$$\begin{array}{|c|c|c|c|}
 \hline
 {\bf 300}&\text{Lower approximation}&\text{Upper approximation}&\text{Upper error}\\
 \hline
 1&17+\frac{11}{35}&17+\frac{11}{34}&\frac{121}{1156}\\
2&17+\frac{22}{69}&17+\frac{11}{34}&\frac{121}{1156}\\
3&17+\frac{33}{103}&17+\frac{11}{34}&\frac{121}{1156}\\
4&17+\frac{33}{103}&17+\frac{44}{137}&\frac{429}{18769}\\
5&17+\frac{33}{103}&17+\frac{77}{240}&\frac{649}{57600}\\
6&17+\frac{33}{103}&17+\frac{110}{343}&\frac{781}{117649}\\
7&17+\frac{33}{103}&17+\frac{143}{446}&\frac{825}{198916}\\
8&17+\frac{33}{103}&17+\frac{176}{549}&\frac{781}{301401}\\
9&17+\frac{33}{103}&17+\frac{209}{652}&\frac{649}{425104}\\
10&17+\frac{33}{103}&17+\frac{242}{755}&\frac{429}{570025}\\
11&17+\frac{33}{103}&{\color{red}17+\frac{25}{78}}&\frac{1}{6084}\\
 \hline
\end{array}$$

$$\begin{array}{|c|c|c|c|}
 \hline
 {\bf 375}&\text{Lower approximation}&\text{Upper approximation}&\text{Upper error}\\
 \hline
1&19+\frac{14}{39}&19+\frac{7}{19}&\frac{49}{361}\\
2&19+\frac{21}{58}&19+\frac{7}{19}&\frac{49}{361}\\
3&19+\frac{4}{11}&19+\frac{7}{19}&\frac{49}{361}\\
4&19+\frac{4}{11}&19+\frac{11}{30}&\frac{61}{900}\\
5&19+\frac{4}{11}&19+\frac{15}{41}&\frac{61}{1681}\\
6&19+\frac{4}{11}&19+\frac{19}{52}&\frac{49}{2704}\\
7&19+\frac{4}{11}&19+\frac{23}{63}&\frac{25}{3969}\\
8&19+\frac{27}{74}&19+\frac{23}{63}&\frac{25}{3969}\\
9&19+\frac{27}{74}&19+\frac{50}{137}&\frac{34}{18769}\\
10&19+\frac{27}{74}&19+\frac{77}{211}&\frac{21}{44521}\\
11&19+\frac{104}{285}&19+\frac{77}{211}&\frac{21}{44521}\\
12&19+\frac{104}{285}&19+\frac{181}{496}&\frac{25}{246016}\\
13&19+\frac{104}{285}&{\color{red}19+\frac{285}{781}}&\frac{1}{609961}\\
 \hline
\end{array}$$

$$\begin{array}{|c|c|c|c|}
 \hline
 {\bf 135}&\text{Lower approximation}&\text{Upper approximation}&\text{Upper error}\\
 \hline
 1&11+\frac{14}{23}&11+\frac{7}{11}&\frac{49}{121}\\
2&11+\frac{21}{34}&11+\frac{7}{11}&\frac{49}{121}\\
3&11+\frac{21}{34}&11+\frac{28}{45}&\frac{154}{2025}\\
4&11+\frac{21}{34}&11+\frac{49}{79}&\frac{189}{6241}\\
5&11+\frac{21}{34}&11+\frac{70}{113}&\frac{154}{12769}\\
6&11+\frac{21}{34}&{\color{red}11+\frac{13}{21}}&\frac{1}{441}\\
 \hline
\end{array}$$

$$\begin{array}{|c|c|c|c|}
 \hline
 {\bf 75}&\text{Lower approximation}&\text{Upper approximation}&\text{Upper error}\\
 \hline
 1&8+\frac{11}{17}&8+\frac{11}{16}&\frac{121}{256}\\
2&8+\frac{11}{17}&8+\frac{2}{3}&\frac{1}{9}\\
3&8+\frac{13}{20}&8+\frac{2}{3}&\frac{1}{9}\\
4&8+\frac{15}{23}&8+\frac{2}{3}&\frac{1}{9}\\
5&8+\frac{17}{26}&8+\frac{2}{3}&\frac{1}{9}\\
6&8+\frac{19}{29}&8+\frac{2}{3}&\frac{1}{9}\\
7&8+\frac{21}{32}&8+\frac{2}{3}&\frac{1}{9}\\
8&8+\frac{23}{35}&8+\frac{2}{3}&\frac{1}{9}\\
9&8+\frac{25}{38}&8+\frac{2}{3}&\frac{1}{9}\\
10&8+\frac{27}{41}&8+\frac{2}{3}&\frac{1}{9}\\
11&8+\frac{29}{44}&8+\frac{2}{3}&\frac{1}{9}\\
12&8+\frac{31}{47}&8+\frac{2}{3}&\frac{1}{9}\\
13&8+\frac{33}{50}&8+\frac{2}{3}&\frac{1}{9}\\
14&8+\frac{33}{50}&8+\frac{35}{53}&\frac{6}{2809}\\
15&8+\frac{68}{103}&8+\frac{35}{53}&\frac{6}{2809}\\
16&8+\frac{68}{103}&{\color{red}8+\frac{103}{156}}&\frac{1}{24336}\\
 \hline
\end{array}$$

$$\begin{array}{|c|c|c|c|}
 \hline
 {\bf 756}&\text{Lower approximation}&\text{Upper approximation}&\text{Upper error}\\
 \hline
 1&27+\frac{27}{55}&27+\frac{1}{2}&\frac{1}{4}\\
2&27+\frac{28}{57}&27+\frac{1}{2}&\frac{1}{4}\\
3&27+\frac{29}{59}&27+\frac{1}{2}&\frac{1}{4}\\
4&27+\frac{30}{61}&27+\frac{1}{2}&\frac{1}{4}\\
5&27+\frac{31}{63}&27+\frac{1}{2}&\frac{1}{4}\\
6&27+\frac{32}{65}&27+\frac{1}{2}&\frac{1}{4}\\
7&27+\frac{33}{67}&27+\frac{1}{2}&\frac{1}{4}\\
8&27+\frac{34}{69}&27+\frac{1}{2}&\frac{1}{4}\\
9&27+\frac{35}{71}&27+\frac{1}{2}&\frac{1}{4}\\
10&27+\frac{36}{73}&27+\frac{1}{2}&\frac{1}{4}\\
11&27+\frac{37}{75}&27+\frac{1}{2}&\frac{1}{4}\\
12&27+\frac{38}{77}&27+\frac{1}{2}&\frac{1}{4}\\
13&27+\frac{39}{79}&27+\frac{1}{2}&\frac{1}{4}\\
14&27+\frac{40}{81}&27+\frac{1}{2}&\frac{1}{4}\\
15&27+\frac{41}{83}&27+\frac{1}{2}&\frac{1}{4}\\
16&27+\frac{42}{85}&27+\frac{1}{2}&\frac{1}{4}\\
17&27+\frac{43}{87}&27+\frac{1}{2}&\frac{1}{4}\\
18&27+\frac{44}{89}&27+\frac{1}{2}&\frac{1}{4}\\
19&27+\frac{45}{91}&27+\frac{1}{2}&\frac{1}{4}\\
20&27+\frac{46}{93}&27+\frac{1}{2}&\frac{1}{4}\\
 \hline
\end{array}$$

$$\begin{array}{|c|c|c|c|}
 \hline
 {\bf 756}&\text{Lower approximation}&\text{Upper approximation}&\text{Upper error}\\
 \hline
21&27+\frac{47}{95}&27+\frac{1}{2}&\frac{1}{4}\\
22&27+\frac{48}{97}&27+\frac{1}{2}&\frac{1}{4}\\
23&27+\frac{49}{99}&27+\frac{1}{2}&\frac{1}{4}\\
24&27+\frac{50}{101}&27+\frac{1}{2}&\frac{1}{4}\\
25&27+\frac{51}{103}&27+\frac{1}{2}&\frac{1}{4}\\
26&27+\frac{52}{105}&27+\frac{1}{2}&\frac{1}{4}\\
27&27+\frac{53}{107}&27+\frac{1}{2}&\frac{1}{4}\\
28&27+\frac{54}{109}&27+\frac{1}{2}&\frac{1}{4}\\
29&27+\frac{54}{109}&27+\frac{55}{111}&\frac{28}{12321}\\
30&27+\frac{54}{109}&{\color{red}27+\frac{109}{220}}&\frac{1}{48400}\\
 \hline
\end{array}$$

$$\begin{array}{|c|c|c|c|}
 \hline
 {\bf 611}&\text{Lower approximation}&\text{Upper approximation}&\text{Upper error}\\
 \hline
 1&24+\frac{5}{7}&24+\frac{35}{48}&\frac{1225}{2304}\\
2&24+\frac{5}{7}&24+\frac{8}{11}&\frac{53}{121}\\
3&24+\frac{5}{7}&24+\frac{13}{18}&\frac{61}{324}\\
4&24+\frac{5}{7}&24+\frac{18}{25}&\frac{49}{625}\\
5&24+\frac{5}{7}&24+\frac{23}{32}&\frac{17}{1024}\\
6&24+\frac{28}{39}&24+\frac{23}{32}&\frac{17}{1024}\\
7&24+\frac{51}{71}&24+\frac{23}{32}&\frac{17}{1024}\\
8&24+\frac{51}{71}&24+\frac{74}{103}&\frac{17}{10609}\\
9&24+\frac{125}{174}&24+\frac{74}{103}&\frac{17}{10609}\\
10&24+\frac{199}{277}&24+\frac{74}{103}&\frac{17}{10609}\\
11&24+\frac{199}{277}&24+\frac{273}{380}&\frac{49}{144400}\\
12&24+\frac{199}{277}&24+\frac{472}{657}&\frac{61}{431649}\\
13&24+\frac{199}{277}&24+\frac{671}{934}&\frac{53}{872356}\\
14&24+\frac{199}{277}&24+\frac{870}{1211}&\frac{25}{1466521}\\
15&24+\frac{1069}{1488}&24+\frac{870}{1211}&\frac{25}{1466521}\\
16&24+\frac{1069}{1488}&24+\frac{1939}{2699}&\frac{14}{7284601}\\
17&24+\frac{3008}{4187}&24+\frac{1939}{2699}&\frac{14}{7284601}\\
18&24+\frac{4947}{6886}&24+\frac{1939}{2699}&\frac{14}{7284601}\\
19&24+\frac{4947}{6886}&{\color{red}24+\frac{6886}{9585}}&\frac{1}{91872225}\\
 \hline
\end{array}$$

$$\begin{array}{|c|c|c|c|}
 \hline
 {\bf 231}&\text{Lower approximation}&\text{Upper approximation}&\text{Upper error}\\
 \hline
 1&15+\frac{6}{31}&15+\frac{1}{5}&\frac{1}{25}\\
2&15+\frac{7}{36}&15+\frac{1}{5}&\frac{1}{25}\\
3&15+\frac{8}{41}&15+\frac{1}{5}&\frac{1}{25}\\
4&15+\frac{9}{46}&15+\frac{1}{5}&\frac{1}{25}\\
5&15+\frac{10}{51}&15+\frac{1}{5}&\frac{1}{25}\\
6&15+\frac{11}{56}&15+\frac{1}{5}&\frac{1}{25}\\
7&15+\frac{12}{61}&15+\frac{1}{5}&\frac{1}{25}\\
8&15+\frac{13}{66}&15+\frac{1}{5}&\frac{1}{25}\\
9&15+\frac{14}{71}&15+\frac{1}{5}&\frac{1}{25}\\
10&15+\frac{15}{76}&15+\frac{1}{5}&\frac{1}{25}\\
 \hline
\end{array}$$

$$\begin{array}{|c|c|c|c|}
 \hline
 {\bf 231}&\text{Lower approximation}&\text{Upper approximation}&\text{Upper error}\\
 \hline
11&15+\frac{16}{81}&15+\frac{1}{5}&\frac{1}{25}\\
12&15+\frac{17}{86}&15+\frac{1}{5}&\frac{1}{25}\\
13&15+\frac{18}{91}&15+\frac{1}{5}&\frac{1}{25}\\
14&15+\frac{19}{96}&15+\frac{1}{5}&\frac{1}{25}\\
15&15+\frac{20}{101}&15+\frac{1}{5}&\frac{1}{25}\\
16&15+\frac{21}{106}&15+\frac{1}{5}&\frac{1}{25}\\
17&15+\frac{22}{111}&15+\frac{1}{5}&\frac{1}{25}\\
18&15+\frac{23}{116}&15+\frac{1}{5}&\frac{1}{25}\\
19&15+\frac{24}{121}&15+\frac{1}{5}&\frac{1}{25}\\
20&15+\frac{25}{126}&15+\frac{1}{5}&\frac{1}{25}\\
21&15+\frac{26}{131}&15+\frac{1}{5}&\frac{1}{25}\\
22&15+\frac{27}{136}&15+\frac{1}{5}&\frac{1}{25}\\
23&15+\frac{28}{141}&15+\frac{1}{5}&\frac{1}{25}\\
24&15+\frac{29}{146}&15+\frac{1}{5}&\frac{1}{25}\\
25&15+\frac{30}{151}&15+\frac{1}{5}&\frac{1}{25}\\
26&15+\frac{30}{151}&15+\frac{31}{156}&\frac{25}{24336}\\
27&15+\frac{30}{151}&15+\frac{61}{307}&\frac{37}{94249}\\
28&15+\frac{30}{151}&15+\frac{91}{458}&\frac{37}{209764}\\
29&15+\frac{30}{151}&15+\frac{121}{609}&\frac{25}{370881}\\
30&15+\frac{30}{151}&{\color{red}15+\frac{151}{760}}&\frac{1}{577600}\\
 \hline
\end{array}$$

 $$\begin{array}{|c|c|c|c|}
 \hline
 {\bf 800}&\text{Lower approximation}&\text{Upper approximation}&\text{Upper error}\\
 \hline
 1&28+\frac{16}{57}&28+\frac{2}{7}&\frac{4}{49}\\
2&28+\frac{9}{32}&28+\frac{2}{7}&\frac{4}{49}\\
3&28+\frac{11}{39}&28+\frac{2}{7}&\frac{4}{49}\\
4&28+\frac{13}{46}&28+\frac{2}{7}&\frac{4}{49}\\
5&28+\frac{15}{53}&28+\frac{2}{7}&\frac{4}{49}\\
6&28+\frac{17}{60}&28+\frac{2}{7}&\frac{4}{49}\\
7&28+\frac{19}{67}&28+\frac{2}{7}&\frac{4}{49}\\
8&28+\frac{21}{74}&28+\frac{2}{7}&\frac{4}{49}\\
9&28+\frac{23}{81}&28+\frac{2}{7}&\frac{4}{49}\\
10&28+\frac{25}{88}&28+\frac{2}{7}&\frac{4}{49}\\
11&28+\frac{27}{95}&28+\frac{2}{7}&\frac{4}{49}\\
12&28+\frac{27}{95}&28+\frac{29}{102}&\frac{25}{10404}\\
13&28+\frac{56}{197}&28+\frac{29}{102}&\frac{25}{10404}\\
14&28+\frac{56}{197}&28+\frac{85}{299}&\frac{49}{89401}\\
15&28+\frac{56}{197}&28+\frac{141}{496}&\frac{41}{246016}\\
16&28+\frac{56}{197}&{\color{red}28+\frac{197}{693}}&\frac{1}{480249}\\
 \hline
\end{array}$$

 $$\begin{array}{|c|c|c|c|}
 \hline
 {\bf 4100}&\text{Lower approximation}&\text{Upper approximation}&\text{Upper error}\\
 \hline
 1&64+\frac{4}{129}&{\color{red}64+\frac{1}{32}}&\frac{1}{1024}\\
 \hline
\end{array}$$

$$\begin{array}{|c|c|c|c|}
 \hline
 {\bf 2000}&\text{Lower approximation}&\text{Upper approximation}&\text{Upper error}\\
 \hline
 1&44+\frac{64}{89}&44+\frac{8}{11}&\frac{64}{121}\\
2&44+\frac{18}{25}&44+\frac{8}{11}&\frac{64}{121}\\
3&44+\frac{18}{25}&44+\frac{13}{18}&\frac{25}{324}\\
4&44+\frac{31}{43}&44+\frac{13}{18}&\frac{25}{324}\\
5&44+\frac{44}{61}&44+\frac{13}{18}&\frac{25}{324}\\
6&44+\frac{44}{61}&44+\frac{57}{79}&\frac{89}{6241}\\
7&44+\frac{44}{61}&44+\frac{101}{140}&\frac{121}{19600}\\
8&44+\frac{44}{61}&44+\frac{145}{201}&\frac{121}{40401}\\
9&44+\frac{44}{61}&{\color{blue}44+\frac{189}{262}}&\frac{89}{68644}\\
10&44+\frac{44}{61}&44+\frac{233}{323}&\frac{25}{104329}\\
11&44+\frac{277}{384}&44+\frac{233}{323}&\frac{25}{104329}\\
12&44+\frac{510}{707}&44+\frac{233}{323}&\frac{25}{104329}\\
13&44+\frac{743}{1030}&44+\frac{233}{323}&\frac{25}{104329}\\
14&44+\frac{743}{1030}&44+\frac{976}{1353}&\frac{64}{1830609}\\
15&44+\frac{743}{1030}&44+\frac{1719}{2383}&\frac{41}{5678689}\\
16&44+\frac{2462}{3413}&44+\frac{1719}{2383}&\frac{41}{5678689}\\
17&44+\frac{2462}{3413}&44+\frac{4181}{5796}&\frac{25}{33593616}\\
18&44+\frac{6643}{9209}&44+\frac{4181}{5796}&\frac{25}{33593616}\\
19&44+\frac{10824}{15005}&44+\frac{4181}{5796}&\frac{25}{33593616}\\
20&44+\frac{10824}{15005}&44+\frac{15005}{20801}&\frac{1}{432681601}\\
 \hline
\end{array}$$

   In the ninth iteration appears $44+\frac{189}{262}$ which is the result after simplify by 11 the approximation $44\frac{2079}{2882}$. The first upper optimal approximation does not appear until the 20th iteration.

$$\begin{array}{|c|c|c|c|}
 \hline
 {\bf 9600}&\text{Lower approximation}&\text{Upper approximation}&\text{Upper error}\\
 \hline
 1&97+\frac{191}{195}&{\color{red}97+\frac{191}{194}}&\frac{36481}{37636}\\
2&97+\frac{191}{195}&97+\frac{382}{389}&\frac{71625}{151321}\\
3&97+\frac{191}{195}&97+\frac{573}{584}&\frac{105241}{341056}\\
4&97+\frac{191}{195}&97+\frac{764}{779}&\frac{137329}{606841}\\
5&97+\frac{191}{195}&97+\frac{955}{974}&\frac{167889}{948676}\\
6&97+\frac{191}{195}&97+\frac{1146}{1169}&\frac{196921}{1366561}\\
7&97+\frac{191}{195}&97+\frac{1337}{1364}&\frac{224425}{1860496}\\
8&97+\frac{191}{195}&97+\frac{1528}{1559}&\frac{250401}{2430481}\\
9&97+\frac{191}{195}&97+\frac{1719}{1754}&\frac{274849}{3076516}\\
10&97+\frac{191}{195}&97+\frac{1910}{1949}&\frac{297769}{3798601}\\
11&97+\frac{191}{195}&97+\frac{2101}{2144}&\frac{319161}{4596736}\\
12&97+\frac{191}{195}&97+\frac{2292}{2339}&\frac{339025}{5470921}\\
13&97+\frac{191}{195}&97+\frac{2483}{2534}&\frac{357361}{6421156}\\
14&97+\frac{191}{195}&97+\frac{2674}{2729}&\frac{374169}{7447441}\\
15&97+\frac{191}{195}&97+\frac{2865}{2924}&\frac{389449}{8549776}\\
16&97+\frac{191}{195}&97+\frac{3056}{3119}&\frac{403201}{9728161}\\
17&97+\frac{191}{195}&97+\frac{3247}{3314}&\frac{415425}{10982596}\\
18&97+\frac{191}{195}&97+\frac{3438}{3509}&\frac{426121}{12313081}\\
19&97+\frac{191}{195}&97+\frac{3629}{3704}&\frac{435289}{13719616}\\
20&97+\frac{191}{195}&97+\frac{3820}{3899}&\frac{442929}{15202201}\\
\hline
\end{array}$$

$$\begin{array}{|c|c|c|c|}
 \hline
 {\bf 9600}&\text{Lower approximation}&\text{Upper approximation}&\text{Upper error}\\
 \hline
21&97+\frac{191}{195}&97+\frac{4011}{4094}&\frac{449041}{16760836}\\
22&97+\frac{191}{195}&97+\frac{4202}{4289}&\frac{453625}{18395521}\\
23&97+\frac{191}{195}&97+\frac{4393}{4484}&\frac{456681}{20106256}\\
24&97+\frac{191}{195}&97+\frac{4584}{4679}&\frac{458209}{21893041}\\
25&97+\frac{191}{195}&97+\frac{4775}{4874}&\frac{458209}{23755876}\\
26&97+\frac{191}{195}&97+\frac{4966}{5069}&\frac{456681}{25694761}\\
27&97+\frac{191}{195}&97+\frac{5157}{5264}&\frac{453625}{27709696}\\
28&97+\frac{191}{195}&97+\frac{5348}{5459}&\frac{449041}{29800681}\\
29&97+\frac{191}{195}&97+\frac{5539}{5654}&\frac{442929}{31967716}\\
30&97+\frac{191}{195}&97+\frac{5730}{5849}&\frac{435289}{34210801}\\
31&97+\frac{191}{195}&97+\frac{5921}{6044}&\frac{426121}{36529936}\\
32&97+\frac{191}{195}&97+\frac{6112}{6239}&\frac{415425}{38925121}\\
33&97+\frac{191}{195}&97+\frac{6303}{6434}&\frac{403201}{41396356}\\
34&97+\frac{191}{195}&97+\frac{6494}{6629}&\frac{389449}{43943641}\\
35&97+\frac{191}{195}&97+\frac{6685}{6824}&\frac{374169}{46566976}\\
36&97+\frac{191}{195}&97+\frac{6876}{7019}&\frac{357361}{49266361}\\
37&97+\frac{191}{195}&97+\frac{7067}{7214}&\frac{339025}{52041796}\\
38&97+\frac{191}{195}&97+\frac{7258}{7409}&\frac{319161}{54893281}\\
39&97+\frac{191}{195}&97+\frac{7449}{7604}&\frac{297769}{57820816}\\
40&97+\frac{191}{195}&97+\frac{7640}{7799}&\frac{274849}{60824401}\\
41&97+\frac{191}{195}&97+\frac{7831}{7994}&\frac{250401}{63904036}\\
42&97+\frac{191}{195}&97+\frac{8022}{8189}&\frac{224425}{67059721}\\
43&97+\frac{191}{195}&97+\frac{8213}{8384}&\frac{196921}{70291456}\\
44&97+\frac{191}{195}&97+\frac{8404}{8579}&\frac{167889}{73599241}\\
45&97+\frac{191}{195}&97+\frac{8595}{8774}&\frac{137329}{76983076}\\
46&97+\frac{191}{195}&97+\frac{8786}{8969}&\frac{105241}{80442961}\\
47&97+\frac{191}{195}&97+\frac{8977}{9164}&\frac{71625}{83978896}\\
48&97+\frac{191}{195}&97+\frac{48}{49}&\frac{1}{2401}\\
 \hline
\end{array}$$
 The approximation  $97\frac{191}{194}$ appears in the first iteration and the first upper optimal approximation appears in the iteration number 48.

$$\begin{array}{|c|c|c|c|}
 \hline
 {\bf 127\frac{3}{11}}&\text{Lower approximation}&\text{Upper approximation}&\text{Upper error}\\
 \hline
 1&11&12&\frac{184}{11}\\
2&11+\frac{0}{1}&11+\frac{1}{2}&\frac{219}{44}\\
3&11+\frac{0}{1}&11+\frac{1}{3}&\frac{116}{99}\\
4&11+\frac{1}{4}&11+\frac{1}{3}&\frac{116}{99}\\
5&11+\frac{1}{4}&{\color{red}11+\frac{2}{7}}&\frac{51}{539}\\
6&11+\frac{3}{11}&11+\frac{2}{7}&\frac{51}{539}\\
7&11+\frac{5}{18}&11+\frac{2}{7}&\frac{51}{539}\\
8&11+\frac{7}{25}&11+\frac{2}{7}&\frac{51}{539}\\
9&11+\frac{9}{32}&11+\frac{2}{7}&\frac{51}{539}\\
10&11+\frac{9}{32}&11+\frac{11}{39}&\frac{200}{16731}\\
11&11+\frac{9}{32}&11+\frac{20}{71}&\frac{211}{55451}\\
12&11+\frac{9}{32}&11+\frac{29}{103}&\frac{84}{116699}\\
13&11+\frac{38}{135}&11+\frac{29}{103}&\frac{84}{116699}\\
\hline
\end{array}$$

\vspace*{-0.5cm}
$$\begin{array}{|c|c|c|c|}
\hline
 {\bf 127\frac{3}{11}}&\text{Lower approximation}&\text{Upper approximation}&\text{Upper error}\\
 \hline
14&11+\frac{67}{238}&11+\frac{29}{103}&\frac{84}{116699}\\
15&11+\frac{67}{238}&11+\frac{96}{341}&\frac{9}{116281}\\
16&11+\frac{163}{579}&11+\frac{96}{341}&\frac{9}{116281}\\
17&11+\frac{163}{579}&11+\frac{259}{920}&\frac{51}{9310400}\\
18&11+\frac{422}{1499}&11+\frac{259}{920}&\frac{51}{9310400}\\
19&11+\frac{681}{2419}&11+\frac{259}{920}&\frac{51}{9310400}\\
20&11+\frac{940}{3339}&11+\frac{259}{920}&\frac{51}{9310400}\\
21&11+\frac{1199}{4259}&11+\frac{259}{920}&\frac{51}{9310400}\\
22&11+\frac{1199}{4259}&11+\frac{1458}{5179}&\frac{219}{295042451}\\
23&11+\frac{1199}{4259}&11+\frac{2657}{9438}&\frac{25}{89075844}\\
24&11+\frac{1199}{4259}&11+\frac{3856}{13697}&\frac{219}{2063685899}\\
25&11+\frac{1199}{4259}&11+\frac{5055}{17956}&\frac{51}{3546597296}\\
26&11+\frac{6254}{22215}&11+\frac{5055}{17956}&\frac{51}{3546597296}\\
27&11+\frac{11309}{40171}&11+\frac{5055}{17956}&\frac{51}{3546597296}\\
28&11+\frac{16364}{58127}&11+\frac{5055}{17956}&\frac{51}{3546597296}\\
29&11+\frac{21419}{76083}&11+\frac{5055}{17956}&\frac{51}{3546597296}\\
30&11+\frac{21419}{76083}&11+\frac{26474}{94039}&\frac{9}{8843333521}\\
31&11+\frac{47893}{170122}&11+\frac{26474}{94039}&\frac{9}{8843333521}\\
32&11+\frac{47893}{170122}&11+\frac{74367}{264161}&\frac{84}{767591373131}\\
33&11+\frac{122260}{434283}&11+\frac{74367}{264161}&\frac{84}{767591373131}\\
34&11+\frac{196627}{698444}&11+\frac{74367}{264161}&\frac{84}{767591373131}\\
35&11+\frac{196627}{698444}&11+\frac{270994}{962605}&\frac{211}{10192692246275}\\
36&11+\frac{196627}{698444}&11+\frac{467621}{1661049}&\frac{200}{30349921584411}\\
37&11+\frac{196627}{698444}&11+\frac{664248}{2359493}&\frac{51}{61239279387539}\\
38&11+\frac{860875}{3057937}&11+\frac{664248}{2359493}&\frac{51}{61239279387539}\\
39&11+\frac{1525123}{5417430}&11+\frac{664248}{2359493}&\frac{51}{61239279387539}\\
40&11+\frac{2189371}{7776923}&11+\frac{664248}{2359493}&\frac{51}{61239279387539}\\
41&11+\frac{2853619}{10136416}&11+\frac{664248}{2359493}&\frac{51}{61239279387539}\\
42&11+\frac{2853619}{10136416}&11+\frac{3517867}{12495909}&\frac{116}{1717625159099091}\\
43&11+\frac{6371486}{22632325}&11+\frac{3517867}{12495909}&\frac{116}{1717625159099091}\\
44&11+\frac{6371486}{22632325}&11+\frac{9889353}{35128234}&\frac{219}{13573921063546316}\\
45&11+\frac{6371486}{22632325}&11+\frac{16260839}{57760559}&\frac{184}{36699103935917291}\\
46&11+\frac{6371486}{22632325}&11+\frac{22632325}{80392884}&\frac{1}{6463015797837456}\\
 \hline
\end{array}$$
   Approximation $11\frac{2}{7}$ appears in the fifth iteration and the first upper approximation is in the iteration number 46.
$$\begin{array}{|c|c|c|c|}
 \hline
 {\bf 5 \frac{1}{3}}&\text{Lower approximation}&\text{Upper approximation}&\text{Upper error}\\
 \hline
 1&2&3&\frac{11}{3}\\
2&2+\frac{0}{1}&2+\frac{1}{2}&\frac{11}{12}\\
3&2+\frac{0}{1}&2+\frac{1}{3}&\frac{1}{9}\\
4&2+\frac{1}{4}&2+\frac{1}{3}&\frac{1}{9}\\
5&2+\frac{2}{7}&2+\frac{1}{3}&\frac{1}{9}\\
6&2+\frac{3}{10}&2+\frac{1}{3}&\frac{1}{9}\\
7&2+\frac{4}{13}&2+\frac{1}{3}&\frac{1}{9}\\
8&2+\frac{4}{13}&2+\frac{5}{16}&\frac{11}{768}\\
9&2+\frac{4}{13}&2+\frac{9}{29}&\frac{11}{2523}\\
10&2+\frac{4}{13}&{\color{red}2+\frac{13}{42}}&\frac{1}{1764}\\
 \hline
\end{array}$$

\newpage

To reach ${\color{red}44\frac{2079}{2882}}$ as approximation of  $\sqrt{2000}$ we can proceed as follows:

The Chuquet's approximation of $\sqrt{5}$ is ${\color{green}\frac{1525}{682}}$. We can write so that $2000=\frac{10000}{5}$, and using the above approximation of $\sqrt{5}$, we get $\frac{100\cdot 682}{1525}=44+\frac{44}{61}$, which is a lower approximation.

$44+\frac{44}{60}=44+\frac{11}{15}$ is an upper approximation. From $44+\frac{44}{61}$ and $44+\frac{11}{15}$, using the rule above exposed and simplifying by 2, we can obtain  $44\frac{2079}{2882}$ in the 48th iteration. The first upper optimal approximation appears in the iteration number 81.

$$\begin{array}{|c|c|c|c|}
 \hline
 {\bf 2000}&\text{Lower approximation}&\text{Upper approximation}&\text{Upper error}\\
 \hline
1&44+\frac{44}{61}&44+\frac{11}{15}&\frac{241}{225}\\
2&44+\frac{44}{61}&44+\frac{55}{76}&\frac{1201}{5776}\\
3&44+\frac{44}{61}&44+\frac{99}{137}&\frac{2129}{18769}\\
4&44+\frac{44}{61}&44+\frac{143}{198}&\frac{25}{324}\\
5&44+\frac{44}{61}&44+\frac{187}{259}&\frac{3889}{67081}\\
6&44+\frac{44}{61}&44+\frac{231}{320}&\frac{4721}{102400}\\
7&44+\frac{44}{61}&44+\frac{275}{381}&\frac{5521}{145161}\\
8&44+\frac{44}{61}&44+\frac{319}{442}&\frac{6289}{195364}\\
9&44+\frac{44}{61}&44+\frac{363}{503}&\frac{7025}{253009}\\
10&44+\frac{44}{61}&44+\frac{407}{564}&\frac{7729}{318096}\\
11&44+\frac{44}{61}&44+\frac{451}{625}&\frac{8401}{390625}\\
12&44+\frac{44}{61}&44+\frac{495}{686}&\frac{9041}{470596}\\
13&44+\frac{44}{61}&44+\frac{539}{747}&\frac{9649}{558009}\\
14&44+\frac{44}{61}&44+\frac{583}{808}&\frac{10225}{652864}\\
15&44+\frac{44}{61}&44+\frac{627}{869}&\frac{89}{6241}\\
16&44+\frac{44}{61}&44+\frac{671}{930}&\frac{11281}{864900}\\
17&44+\frac{44}{61}&44+\frac{715}{991}&\frac{11761}{982081}\\
18&44+\frac{44}{61}&44+\frac{759}{1052}&\frac{12209}{1106704}\\
19&44+\frac{44}{61}&44+\frac{803}{1113}&\frac{12625}{1238769}\\
20&44+\frac{44}{61}&44+\frac{847}{1174}&\frac{13009}{1378276}\\
21&44+\frac{44}{61}&44+\frac{891}{1235}&\frac{13361}{1525225}\\
22&44+\frac{44}{61}&44+\frac{935}{1296}&\frac{13681}{1679616}\\
23&44+\frac{44}{61}&44+\frac{979}{1357}&\frac{13969}{1841449}\\
24&44+\frac{44}{61}&44+\frac{1023}{1418}&\frac{14225}{2010724}\\
25&44+\frac{44}{61}&44+\frac{1067}{1479}&\frac{14449}{2187441}\\
26&44+\frac{44}{61}&44+\frac{1111}{1540}&\frac{121}{19600}\\
27&44+\frac{44}{61}&44+\frac{1155}{1601}&\frac{14801}{2563201}\\
28&44+\frac{44}{61}&44+\frac{1199}{1662}&\frac{14929}{2762244}\\
29&44+\frac{44}{61}&44+\frac{1243}{1723}&\frac{15025}{2968729}\\
30&44+\frac{44}{61}&44+\frac{1287}{1784}&\frac{15089}{3182656}\\
31&44+\frac{44}{61}&44+\frac{1331}{1845}&\frac{15121}{3404025}\\
32&44+\frac{44}{61}&44+\frac{1375}{1906}&\frac{15121}{3632836}\\
33&44+\frac{44}{61}&44+\frac{1419}{1967}&\frac{15089}{3869089}\\
34&44+\frac{44}{61}&44+\frac{1463}{2028}&\frac{15025}{4112784}\\
35&44+\frac{44}{61}&44+\frac{1507}{2089}&\frac{14929}{4363921}\\
\hline
\end{array}$$

$$\begin{array}{|c|c|c|c|}
\hline
 {\bf 2000}&\text{Lower approximation}&\text{Upper approximation}&\text{Upper error}\\
 \hline
36&44+\frac{44}{61}&44+\frac{1551}{2150}&\frac{14801}{4622500}\\
37&44+\frac{44}{61}&44+\frac{1595}{2211}&\frac{121}{40401}\\
38&44+\frac{44}{61}&44+\frac{1639}{2272}&\frac{14449}{5161984}\\
39&44+\frac{44}{61}&44+\frac{1683}{2333}&\frac{14225}{5442889}\\
40&44+\frac{44}{61}&44+\frac{1727}{2394}&\frac{13969}{5731236}\\
41&44+\frac{44}{61}&44+\frac{1771}{2455}&\frac{13681}{6027025}\\
42&44+\frac{44}{61}&44+\frac{1815}{2516}&\frac{13361}{6330256}\\
43&44+\frac{44}{61}&44+\frac{1859}{2577}&\frac{13009}{6640929}\\
44&44+\frac{44}{61}&44+\frac{1903}{2638}&\frac{12625}{6959044}\\
45&44+\frac{44}{61}&44+\frac{1947}{2699}&\frac{12209}{7284601}\\
46&44+\frac{44}{61}&44+\frac{1991}{2760}&\frac{11761}{7617600}\\
47&44+\frac{44}{61}&44+\frac{2035}{2821}&\frac{11281}{7958041}\\
48&44+\frac{44}{61}&{\color{red}44+\frac{2079}{2882}}&\frac{89}{68644}\\
49&44+\frac{44}{61}&44+\frac{2123}{2943}&\frac{10225}{8661249}\\
50&44+\frac{44}{61}&44+\frac{2167}{3004}&\frac{9649}{9024016}\\
51&44+\frac{44}{61}&44+\frac{2211}{3065}&\frac{9041}{9394225}\\
52&44+\frac{44}{61}&44+\frac{2255}{3126}&\frac{8401}{9771876}\\
53&44+\frac{44}{61}&44+\frac{2299}{3187}&\frac{7729}{10156969}\\
54&44+\frac{44}{61}&44+\frac{2343}{3248}&\frac{7025}{10549504}\\
55&44+\frac{44}{61}&44+\frac{2387}{3309}&\frac{6289}{10949481}\\
56&44+\frac{44}{61}&44+\frac{2431}{3370}&\frac{5521}{11356900}\\
57&44+\frac{44}{61}&44+\frac{2475}{3431}&\frac{4721}{11771761}\\
58&44+\frac{44}{61}&44+\frac{2519}{3492}&\frac{3889}{12194064}\\
59&44+\frac{44}{61}&44+\frac{2563}{3553}&\frac{25}{104329}\\
60&44+\frac{44}{61}&44+\frac{2607}{3614}&\frac{2129}{13060996}\\
61&44+\frac{44}{61}&44+\frac{2651}{3675}&\frac{1201}{13505625}\\
62&44+\frac{44}{61}&44+\frac{2695}{3736}&\frac{241}{13957696}\\
63&44+\frac{2739}{3797}&44+\frac{2695}{3736}&\frac{241}{13957696}\\
64&44+\frac{5434}{7533}&44+\frac{2695}{3736}&\frac{241}{13957696}\\
65&44+\frac{8129}{11269}&44+\frac{2695}{3736}&\frac{241}{13957696}\\
66&44+\frac{10824}{15005}&44+\frac{2695}{3736}&\frac{241}{13957696}\\
67&44+\frac{10824}{15005}&44+\frac{13519}{18741}&\frac{1129}{351225081}\\
68&44+\frac{10824}{15005}&44+\frac{24343}{33746}&\frac{1889}{1138792516}\\
69&44+\frac{10824}{15005}&44+\frac{35167}{48751}&\frac{2521}{2376660001}\\
70&44+\frac{10824}{15005}&44+\frac{45991}{63756}&\frac{25}{33593616}\\
71&44+\frac{10824}{15005}&44+\frac{56815}{78761}&\frac{3401}{6203295121}\\
72&44+\frac{10824}{15005}&44+\frac{67639}{93766}&\frac{3649}{8792062756}\\
73&44+\frac{10824}{15005}&44+\frac{78463}{108771}&\frac{3769}{11831130441}\\
74&44+\frac{10824}{15005}&44+\frac{89287}{123776}&\frac{3761}{15320498176}\\
75&44+\frac{10824}{15005}&44+\frac{100111}{138781}&\frac{3625}{19260165961}\\
76&44+\frac{10824}{15005}&44+\frac{110935}{153786}&\frac{3361}{23650133796}\\
77&44+\frac{10824}{15005}&44+\frac{121759}{168791}&\frac{2969}{28490401681}\\
78&44+\frac{10824}{15005}&44+\frac{132583}{183796}&\frac{2449}{33780969616}\\
79&44+\frac{10824}{15005}&44+\frac{143407}{198801}&\frac{1801}{39521837601}\\
80&44+\frac{10824}{15005}&44+\frac{154231}{213806}&\frac{1025}{45713005636}\\
81&44+\frac{10824}{15005}&44+\frac{165055}{228811}&\frac{1}{432681601}\\
 \hline
\end{array}$$

\section{Conclusions}

It is clear that with this method we can obtain all approximations. They are consistent  with the mathematical knowledge at that time and may be with the Ortega's knowledge, because the book was printed in Lyon, the same place where Chuquet was living and the same city in which he published his text ``Triparty" where  we can find the ``regle des nombres mohines" for computing approximations of square roots.

Nevertheless, computing $\sqrt{2000}$ (see appendix III) is too long and approximates $\sqrt{5\frac 13}$ by $2+\frac 16+\frac 17$ instead of $2\frac{13}{42}$, may suggest some doubt about this hypothesis. In any case, consistency and accuracy of the results leads us to refuse that Ortega or the publisher had made a mistake in the approximations as some authors maintain.

\section{Appendix I. EDITIONS OF THE WORK.}

\begin{itemize}
\item 1512. Lyon

 ``Siguese una composicion de la arte de la aritmetica y Juntamente de geometria: fecha y ordenada por fray Juan de ortega de la orden de santo domingo: de los predicadores."

 Imprimido a Leon : en casa de maistro Nicolau de Benedictis : por Joannes  trinxer librero de barcelona

Reference: Fondo Hist\'orico de la Universidad de Salamanca

Free access digitized copy:

\url{https://gredos.usal.es/handle/10366/83271}

\item 1515 en Lyon. French translation by Claude Plantin.

Auteur(s) :  Ortega, Juan (O.F.P., Le P.)

Oeuvre tres subtille et profitable de l'art de science de aristm\'eticque et g\'eom\'etrie, translat\'e nouvellement d'espaignol en fran\c{c}oys [de fr\`ere Jehan de Lortie, de l'ordre Sainct Dominicque]\ldots Ayez ce livre, n'y faillez nullement ; Symon Vincent si vous en fournira, en rue Merci\`ere o\`u il est demourant\ldots -

``A la fin" : Imprim\'e a Lyon, par maistre Estienne Baland, l'an mil cincq cens et quinze, le XXIII. jour de octobre

 Traduit par fr\`ere Claude Platin, humble religieux de l'ordre de Sainct Anthoine en Viennoys

Reference: Biblioth\`eque nationale de France

\url{https://catalogue.bnf.fr/ark:/12148/cb31041178s.public}

\item 1515.Roma

 Suma de arithmetica, geometria pratica, utilissima, ordinata per Johane de Ortega, Spagnolo Palentino.

Impresso in Roma: per Maestro Stephano Guillieri de Lorena, anno del nostro Signor 1515 adi 10 de Noue[m]bre regnante Leone Papa decimo in suo anno tertio.

    Although this work has a Latin title, it is actually the author's Italian translation and adaptation of his Spanish original.

References:
\begin{itemize}
\item
\url{https://ccuc.csuc.cat/}
\item
\url{https://catalogue.bnf.fr/ark:/12148/cb31041178s.public}
\item
 \url{http://clio.cul.columbia.edu:7018}
\item Free access digitized copy:

\url{https://atena.beic.it/R/}

Searching for: Suma de arithmetica: geometria pratica: vtilissima,
\vspace{.3 cm}

Hugo Marquant in \cite{mar}

\textit{``goes more
concretely (ad intra and extra) into the question of the proper nature of the Italian version (1515) of the Castilian original
(1512) of his famous manual of commercial arithmetic as a kind of 'autotranslation'."}

In the colophon appears

Et perche per li exempli passati poterai fare qual si voglia figura
o mesura che sia simile:non voglio essere piu prolixo:perche si desideri piu exempli li trouerai in vna altra mia opera:laquale composi in, Spagnia:laquale tracta splendidamente:cossi de la Arismetica:como de la Geometria:
per tanto non voglio piu; elargarme \ldots

 \end{itemize}

\item 1522 in Mesina.

Italian edition with remarquable differences from the previous ones. Starting from the title in latin:

 Sequitur la quarta opera de arithmetica \& geometria / facta et ordinata per Johanne de
Ortega ...

and followed by different geometric examples. See table(\ref{em})

Sequitur la quarta opera de arithmetica [et] geometria facta et ordinata per Ioanne de Ortega spagnuolo palentino, la quale fu composta in Messina in lo anno MDXXII regnante lo sanctissimo catholico imperatore d. Carlo re di Spagna vtriusque Sicilie et Ierusalem in lo suo tercio anno in lo tempo de lo summo pontifice papa Adriano sesto. Stampata in la nobili chitati di Misina per Paolo Gottardo Pontio ... a instantia di Pietro Tini 1585

Reference:
Real Biblioteca. Patrimonio Nacional.

\url{http://inventarios.realbiblioteca.es/node/19816}

Exemplars:
\begin{itemize}
\item Real Biblioteca del Monasterio de San Lorenzo de El Escorial

\url{https://rbme--patrimonionacional--es.insuit.net/}

\item Biblioteca Antica del Seminario Vescovile di Padova

\url{http://www.bibliotecaseminariopda.it/}

\item Biblioteca de la Universidad de Columbia

 \url{http://clio.cul.columbia.edu}

\item Biblioteca Mazarino. Par\'is.

\url{http://www.sudoc.abes.fr}
\end{itemize}

\item 1534. Sevilla.

Tratado subtilissimo de Arismetica y de Geometria
c\~opuesto y ordenado por el reuerendo padre fray Juan de Ortega de la orden de los predicadores

En ... Seuilla en casa de Ju\~a Cr\~oberger

Reference: Biblioteca Nacional de Espa\~na. (BNE).

\url{http://datos.bne.es/edicion/a5162736.html}

\item 1537. Sevilla.

Tratado subtilissimo de arismetica y de geometria c\~opuesto y ordenado por el reuer~endo padre fray Juan de Ortega de la orden de los predicadores.

Agora nueuam~ete corregido y emendado En ... Seuilla en casa de Ju\~a Cr\~oberger

References:
\begin{itemize}

\item Biblioteca Nacional de Espa\~na. (BNE).

\url{http://datos.bne.es/edicion/bima0000013858.html}

Free access digitized copy:

\url{http://bdh-rd.bne.es/viewer.vm?id=0000115463&page=1}

\item books.google

Free access digitized copy:

\url{https://books.google.es/books/ucm?id=e5NmUB8tqXYC&printsec=frontcover&hl=es&source=gbs_ge_summary_r&cad=0#v=onepage&q&f=false}

\item books.google

Free access digitized copy:

\url{https://books.google.es/books/about/Tratado_subtilissimo_de_Aritmetica_y_de.html?id=UN9hHQAACAAJ&redir_esc=y}
\end{itemize}

\item 1542. Sevilla.
Tratado subtilissimo de arismetica y de geometria c\~opuesto y ordenado por el reuerendo padre fray Ju\~a de Ortega de la orden de los predicadores.

   Fue impresso el presente libro ... agora nueum\~ete corregido y emendado en casa \~d Jacom Cr\~oberger en la muy noble y muy leal ciudad de Seuilla, 1542.

References;

\begin{itemize}
\item
\url{http://cisne.sim.ucm.es/record=b2338782*spi}
\item
\url{http://clio.cul.columbia.edu:7018}

\item books.google.es

Free access digitized copy:

\url{http://books.google.com/books/ucm?vid=UCM5322482689& printsec=frontcover}
\item
\url{https://ccuc.csuc.cat/}
\end{itemize}

\item 1552. Sevilla.

Tractado subtilissimo d'arismetica y de geometria, compuesto por  el reuer~edo padre fray Juan de Hortega de la orden de los predicadores.
Ahora de nuevo enmendado \ldots por Gon\c{c}alo Busto.

    Fue impresso \~e la muy noble  muy leal ciudad de Seuilla,
     por Ju\~a canalla 

     Acabose \ldots a\~no de nuestro criador y red\~eptor Jesu Christo de
      mill \newline quinientos  cinquenta y dos a\~nos \ldots 1552.

References:

\begin{itemize}
\item
\url{https://ccuc.csuc.cat/}
\item  Biblioteca Nacional de Espa\~na. (BNE).
\item
 \url{http://clio.cul.columbia.edu:7018}
\item
 Biblioteca Nacional de Portugal (BNP).

\item books.google

Free access digitized copy:

\url{https://books.google.es/books/ucm?id=E8XpbKMFDbYC&printsec=frontcover&dq=inauthor:%22Juan+de+Ortega+(O.P.)%22&hl=es&cd=1#v=onepage&q&f=false}
\end{itemize}

\item 1563. Granada

Tractado subtilissimo \~d arismetica y geometria  c\~opuesto por el reuer\~edo padre fray Ju\~a de Hortega ; agora de nueuo emendado \~o mucha dilig\~ecia por Juan Lagarto y antes por Gon\c{c}alo Busto de muchos errores \~q auia en algunas impressiones passadas ; van annadidas en esta impression las prueuas desde reduzir hasta partir quebrados, y en las mas de las figuras de geometria sus prueuas, con ciertos auisos subjectos al algebra. Va a\~nadido en esta postrera impressi\~o vn Tractado del bachiller Iu\~a Perez de Moya : trata reglas para c\~otar sin pluma y de reduzir vnas monedas castellanas en otras

Fue impresso en la muy noble, n\~obrada  gr\~a ciudad de Granada : en casa de Rene Rabut, impressor de libros, junto alos hospitales del Corpus Christi : a costa de Iu\~a dias mercader de libros, 1563, en ocho dias del mes de abril

References:
\begin{itemize}
\item
Fondo Hist\'orico de la Universidad de Salamanca.

\url{http://brumario.usal.es/}
\item
 Universitat de Val\`encia.

 \url{http://trobes.uv.es/search~S1*spi?/aortega+juan/aortega+juan
 /1%2C5%2C6%2CB/frameset&FF=aortega+juan+de+o+p&2%2C%2C2}

 \item
\url{http://clio.cul.columbia.edu:7018/vwebv/holdingsInfo?bibId=1231391}
\end{itemize}

\end{itemize}

\subsection{The hypothetical edition of Cambray of 1612}\

\vspace{0.3cm}

 Zarco del Valle y Sancho Ray\'on  \cite{zar},
 wrote the book

 \vspace{0.3cm}

 \textit{``Ensayo de una biblioteca espa\~nola de libros raros y curiosos"}.

\vspace{0.3cm}

 In this book we find:

\vspace{.3cm}

\textit{ ``ORTEGA (j. DE). Tratado de aritm\'etica, en Cambray, a\~no 1612. (Aa, 184.)"}

\vspace{.3cm}

\begin{figure}[b]
\includegraphics[scale=0.3]{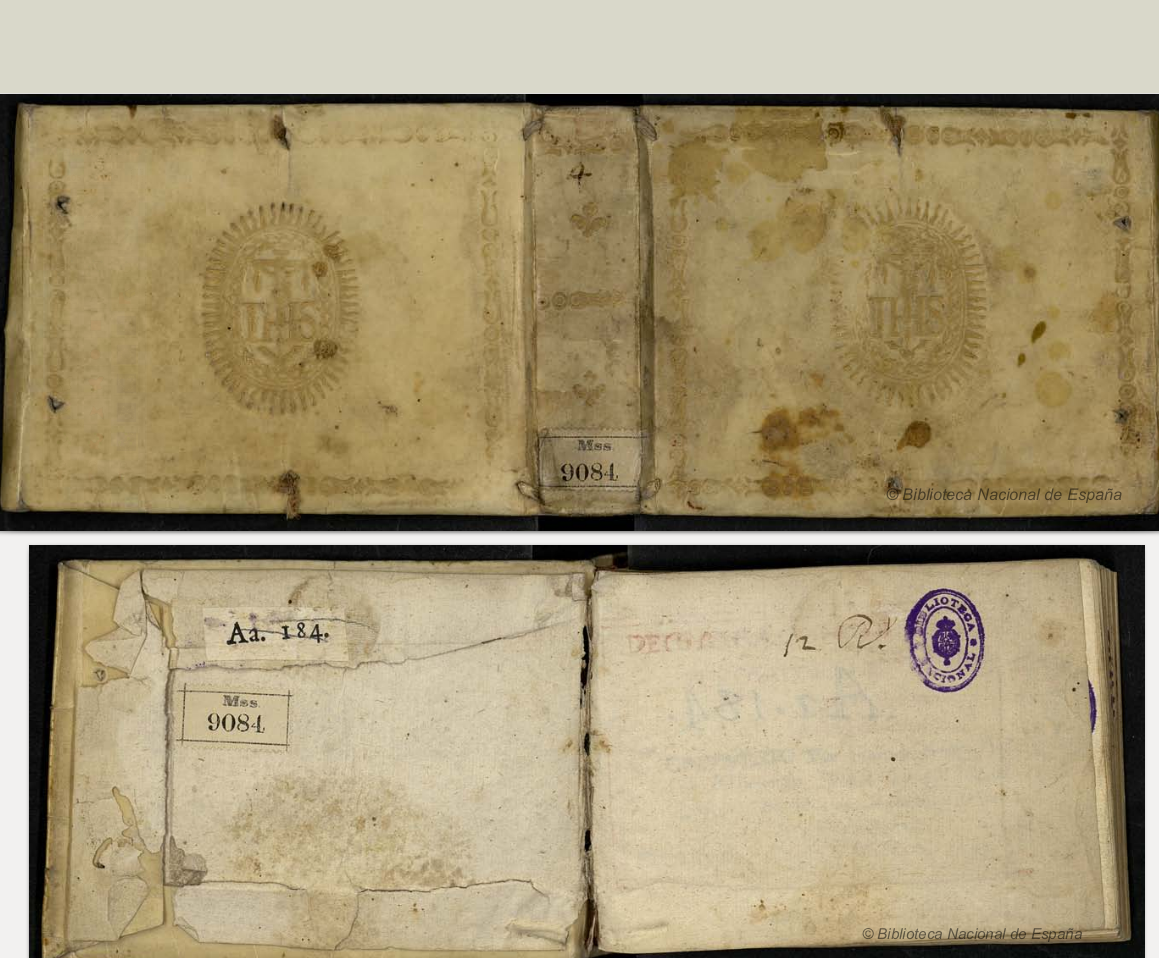}
\caption{Manuscrito 9084.\label{m9084a}}
\end{figure}

This work corresponds to the manuscript Mss 9084 of the  ``Biblioteca Nacional de Espa\~na" (The reference Aa 184 appears in figure \ref{m9084a}).

\vspace{.3cm}
\url{http://bdh.bne.es/bnesearch/CompleteSearch.do?languageView=es&field=autor&text=Juan+de+Ortega&showYearItems=&exact=&textH=&advanced=&completeText=&pageSize=1&pageSizeAbrv=30&pageNumber=18}
\vspace{.3cm}

\begin{figure}
\includegraphics[scale=0.3]{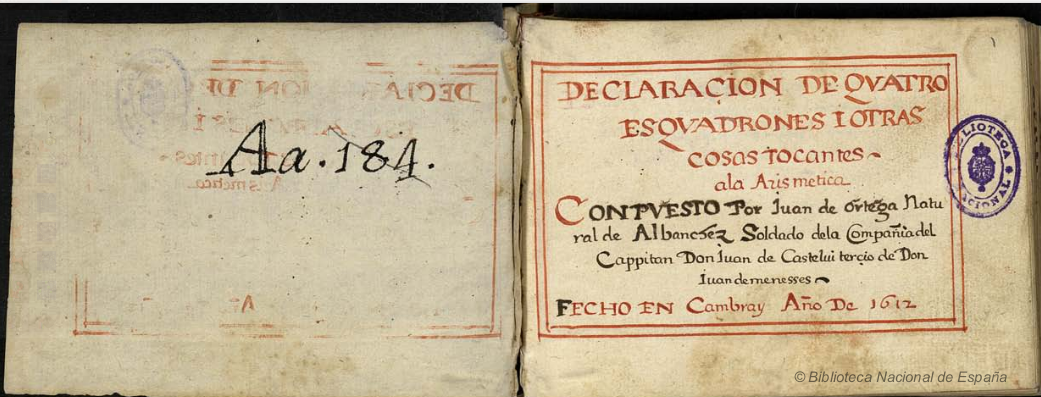}
\caption{Manuscrito 9084.\label{m9084b}}
\end{figure}

\vspace{.3cm}

In the figure \ref{m9084b} we see that its author is \textit {``Juan de Ortega natural of Albanchez"} and its title is \textit{``DECLARACION DE QUATRO ESCUADRONES I OTRAS cosas tocantes a la Arismetica. FECHO EN Cambray A\~no De 1612".}

\newpage


\section{Appendix II.{ Problems in which appear $\sqrt{127\frac{3}{11}}$ and $\sqrt{5\frac 13}$.}}

In the 1512 edition there are two upper approximations that remain unchanged in the following editions. The approximations are: $\sqrt{127\frac{3}{11}}\simeq 11\frac{2}{7}$  and $\sqrt{5\frac 13}\simeq 2+\frac 16+\frac 17$,  we can find them in page 230:

\begin{center}
\includegraphics[scale=0.9]{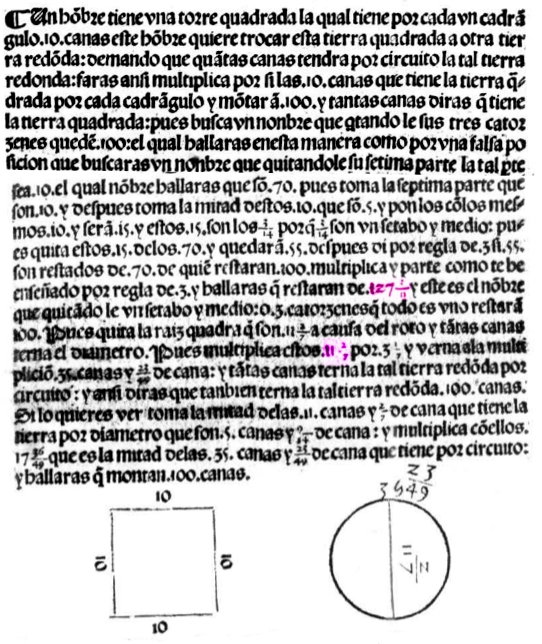}
\end{center}

Un h\~obre tiene vna torre quadrada la qual tiene por cada vn
cadr\~agulo.10.canas este h\~obre quiere trocar esta tierra quadrada a otra tierra red\~oda:demando que qu\~atas canas terna por circuito la tal tierra red\~oda:faras ansi multiplica por si las.10.canas que tiene la tierra \~qdrada por cada cadr\~agulo y m\~otar\~a.100.y tantas canas diras \~q tiene la tierra quadrada:pues busca vn nonbre que qtando le sus tres catorzenes qued\~e.100:el qual hallaras enesta manera como por vna falsa posicion que buscaras vn nonbre que quitandole su septima parte la tal \underline{p}te sea.10.el qual n\~obre hallaras que s\~o.70.pues toma la septima parte que son.10.y despues
toma la mitad destos.10.que s\~o.5.y ponlos c\~olos mesmos.10.y ser\~a.15.y estos.15.son los $\frac{3}{14}$ por\~q $\frac{3}{14}$ son vn setabo y medio:pues quita estos.15.delos.70.y quedar\~a.55.despues di por regla de.3si.55.son restados de.70.de qui\~e restaran.100.multiplica y parte como te he ense\~nado por regla de.3.y hallaras \~q restaran de.{\color{red}$127\frac{3}{11}$}y este es el n\~obre que quit\~adole vn setabo y medio:o.3.catorzenes \~q todo es vno restar\~a.100. Pues quita la raiz quadrada \~q son.{\color{red}$11\frac 27$} a causa del roto y t\~atas canas terna el diametro. Pues multiplica estos.$11\frac 27$por.$3\frac 17$ y verna ala multiplicaci\~o.$35$.canas y $\frac{23}{49}$ de cana:y t\~atas canas terna la tal tierra red\~oda por circuito:y ansi diras que tanbien terna la tal tierra red\~oda.100.canas. Si lo quieres ver toma la mitad delas.11.canas y$\frac 27$de cana que tiene la tierra por diametro que son.5.canas y$\frac {9}{14}$de cara:y multiplica c\~oellos.$17\frac{36}{49}$que es la mitad delas.35.canas y$\frac{23}{49}$de cana que tiene por circuito:y hallaras \~q montan.100.canas.

\vspace{1cm}
\begin{center}
\includegraphics[scale=0.9]{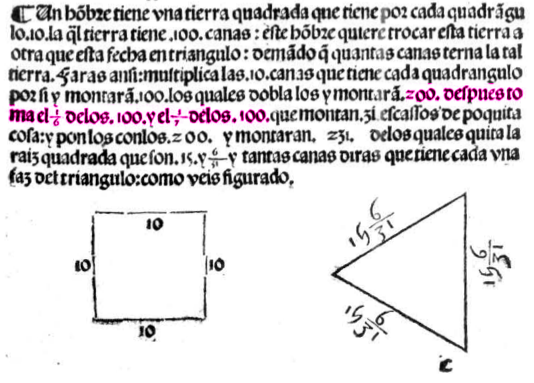}
\end{center}

Un h\~obre tiene vna tierra quadrada que tiene por cada quadr\~agulo.10.la \~ql tierra tiene.100.canas: este h\~obre quiere trocar esta tierra a otra que esta fecha en triangulo:dem\~ado \~q quantas canas terna la tal tierra. Faras ansi:multiplica las.10.canas que tiene cada quadrangulo por si y montar\~a.100.los quales dobla los y montar\~a.{\color{red}$200$.despues toma el$\frac{1}{6}$delos.100.y el$\frac 17$ delos.100.}que montan.31 escassos de poquita cosa:y ponlos conlos.200. y montaran. 231. delos quales quita la raiz quadrada que son.$15$.y$\frac{6}{31}$ y tantas canas diras que tiene cada vna faz del triangulo:como veis figurado.
 \begin{center}
\includegraphics[scale=0.5]{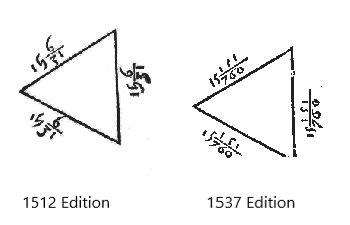}
\end{center}

To compute the square of the edge of an equilateral triangle as a function of his area he has multiply the area by $\frac{4}{\sqrt{3}}=\sqrt{5\frac {1}{3}}$, and to do that he used the approximation $2+\frac{1}{6}+\frac{1}{7}$ in all editions.

\section{APPENDIX III. The approximation of $\sqrt{2000}$}

We can use the result:

$44+\frac{44}{61}<\sqrt{2000}<44+\frac{11}{15}$,
hence, the ``mediation" between

 $44+\frac{44}{61}$ and $44+\frac{11}{15}$ is $44+\frac{55}{76}$

To check if $44+\frac{55}{76}$ is a lower or upper approximation , we can see if  $(44+\frac{55}{76})^2-2000$ is less or grater than 0, and we do it in the following way
$$\frac{55}{76}-\frac{44}{61}=\frac{11}{61\cdot 76}$$
$$(44+\frac{55}{76})^2-2000=(44+\frac{44}{61}+\frac{11}{61\cdot 76})^2-2000=-\frac{16}{61^2}+\frac{5456\cdot 11}{61^2}\cdot \frac{1}{76}+\frac{11^2}{61^2}\cdot \frac{1}{76^2}>0$$
since $\frac{5456\cdot 11}{76}>16$

Since $\frac{5456\cdot 11}{137}>16$ and the ``mediation" between  $44+\frac{44}{61}$ and
$44+\frac{55}{76}$ is $44+\frac{99}{137}$.

$$ \frac{99}{137}-\frac{44}{61}=\frac{11}{61\cdot 137}$$

$$(44+\frac{99}{137})^2-2000=(44+\frac{44}{61}+\frac{11}{61\cdot 137})^2-2000=$$
$$=-\frac{16}{61^2}+\frac{5456\cdot 11}{61^2}\cdot \frac{1}{137}+\frac{11^2}{61^2}\cdot \frac{1}{137^2}>0$$

$$(44+\frac{44}{61}+\frac{11}{61\cdot x})^2-2000=-\frac{16}{61^2}+\frac{5456\cdot 11}{61^2}\cdot \frac{1}{x}+\frac{11^2}{61^2}\cdot \frac{1}{x^2}$$

And if $x<\frac{5456\cdot 11}{16}= 3751$ the approximation $44+\frac{44}{61}+\frac{11}{61\cdot x}$ is an upper approximation of $\sqrt{2000}$.

So that, all ``mediations" to reach
$44+\frac{2079}{2882}$ are upper aproximations.


\end{document}